\documentclass[journal]{IEEEtran}
\usepackage[pdftex]{graphicx} 
\usepackage{amsmath} 
\usepackage{algorithm,algorithmic} 
\usepackage[caption=false,font=footnotesize]{subfig} 
\usepackage{url} 
\usepackage{xcolor}
\def\ones { \mathbf{1} }

\def\IR { \textrm{I}\!\textrm{R}} %
\def\IP { \textrm{I}\!\textrm{P}} %

\def\A { \mathbf{A}}

\def\D { \mathbf{D}}

\def\H { \mathbf{H}} 
\def\I { \mathbf{I}}

\def\N { \mathbf{N}}

\def\Q { \mathbf{Q}}

\def\U { \mathbf{U}}
\def\V { \mathbf{V}}
\def\W { \mathbf{W}}
\def\X { \mathbf{X}} 
\def\Y { \mathbf{Y}}

\def\h { \mathbf{h}}

\def\p { \mathbf{p}}

\def\w { \mathbf{w}}
\def\x { \mathbf{x}} 
\def\y { \mathbf{y}}

\begin{document} 
\title{Algorithms and Comparisons of Nonnegative Matrix Factorizations with Volume Regularization for Hyperspectral Unmixing} 
\author{Andersen~Man~Shun~Ang,~\IEEEmembership{Member,~IEEE,}
		and~Nicolas~Gillis,~\IEEEmembership{Member,~IEEE}
		\thanks{A. Ang and N. Gillis are with the Department of Mathematics and Operational Research,
				Universit\'e de Mons, Belgium. E-mails: \{manshun.ang,nicolas.gillis\}@umons.ac.be. 
				This work was supported by the European Research Council (ERC starting grant no 679515), and 	
				the Fonds de la Recherche Scientifique - FNRS and the Fonds Wetenschappelijk Onderzoek - Vlanderen (FWO) under EOS Project no O005318F-RG47.}
	    }
\maketitle
\begin{abstract} 
In this work, we consider nonnegative matrix factorization (NMF) with a regularization that promotes small volume of the convex hull spanned by the basis matrix. 
We present highly efficient algorithms for three different volume regularizers, and compare them on endmember recovery in hyperspectral unmixing. 
The NMF algorithms developed in this work are shown to outperform the state-of-the-art volume-regularized NMF methods, and produce meaningful decompositions on real-world hyperspectral images in situations where endmembers are highly mixed (no pure pixels).  
Furthermore, our extensive numerical experiments show that when the data is highly separable, meaning that there are data points close to the true endmembers, and there are a few endmembers, the regularizer based on the determinant of the Gramian produces the best results in most cases. 
For data that is less separable and/or contains more endmembers, the regularizer based on the logarithm of the determinant of the Gramian performs best in general. 
\end{abstract}

	

\section{Introduction} \label{intro}

\IEEEPARstart{N}on-negative Matrix Factorization (NMF) is the following problem: given a matrix $\X \in \IR^{m \times n}$ and an integer $r$, 
find two matrices $\W \in \IR^{m \times r}_+$ and $\H \in \IR^{r \times n}_+$ such that 
\vspace*{-1mm}
\begin{equation} 
\label{NMF}
(\textrm{NMF}) : ~ \X \cong \W\H. 
\vspace*{-1mm}
\end{equation}
NMF has many applications; see, e.g.,~\cite{cichocki2009nonnegative, gillis2014and, fu2018nonnegative} and the references therein. 
In this work, we focus on blind {hyperspectral unmixing} (HU)~\cite{bioucas2013hyperspectral} that aims at recovering the spectral signatures of the pure materials (called endmembers, represented as the columns of $\W$) and their abundances in each pixel (represented as the columns of $\H$) within a hyperspectral image $\X$ where each column of $\X$ is the spectral signature of a pixel; see Section~\ref{Intro_background} for more details. 
In order to tackle HU, we will use \textit{volume-regularized NMF} (VRNMF) which can be formulated as follows 
\vspace*{-1mm}
\begin{eqnarray}
\min_{\W,\H}  & & f(\W,\H; \X) + \lambda V(\W)   \nonumber \\
\textrm{subject to } & & \W \geq 0, ~ \H \geq 0, ~ \H^\top\ones_r \leq \ones_n, \label{VRNMF_pbl}  
\vspace*{-1mm}
\end{eqnarray}
where $\ones_n$ denotes the vector of ones of length $n$ and $\W \geq 0$ indicates that $\W$ is component-wise non-negative.
The regularization parameter $\lambda \geq 0$ controls the trade off between the data fitting term $f(\W,\H; \X)$ and the volume regularizer $V(\W)$. 
The constraint $\H^\top\ones_r \leq \ones_n$ is a relaxation of the sum-to-one constraint $\H^\top\ones_r = \ones_n$ which requires that the abundances of the endmembers in each pixel sum to one. 
Because VRNMF minimizes the volume of $\W$, the constraint $\H^\top\ones_r \leq \ones_n$ will tend to be active for most pixels.  
It has been shown previously that this relaxation works better in practice as it allows to take into account for example different illumination conditions within the hyperspectral image; see for example~\cite{gillis2014successive}. The reason to consider VRNMF has a long history in blind HU and is motivated by geometric insights: the columns of $\W$ (the endmembers) are the vertices of a convex hull that contains the data points. In the absence of pure pixels, that is, pixels containing a single endmember, minimizing the volume of the columns of $\W$ allows to recover these endmemebers under mild conditions; see Figure~\ref{Fig_Toy} for an illustation, and Section~\ref{Intro_background} for more details. 

In this paper, we will consider the most widely used data fitting term, namely the least squares error  
$f(\W,\H; \X) 
= \frac{1}{2} \| \X-\W\H\|_F^2 
=  \frac{1}{2} \sum_{i,j} (\X-\W\H)_{ij}^2$. 
For the volume regularizer, we will consider the following three functions 
\begin{eqnarray}
\text{Detminant:}  & V_\text{det}(\W) & = \frac{1}{2}\det (\W^\top\W), \nonumber\\
\text{Log-det:}    & V_\text{logdet}(\W) & = \frac{1}{2}\log\det (\W^\top\W + \delta \I_r), \nonumber \\
\text{Nuclear:}    & V_*(\W) & = \| \W \|_*,   \nonumber
\end{eqnarray}
where $\det(\A)$ is the determinant of matrix $\A$, $\I_r$ is identity matrix of size $r$, $\delta > 0$ is a constant to lower bound $V_\text{logdet}$, and $\|\A\|_*$ is the nuclear norm of $\A$, that is, the sum of the singular values of $\A$. 
VRNMF aims at fitting the data points within the convex hull of the columns of $\W$ which should have a small volume. 
The reason to consider the functions $V_\text{det}$ and $V_\text{logdet}$ is that $\sqrt{\text{det}(\W^\top\W)}/r!$ is the volume of the convex hull of the columns of $\W$ and the origin. 
Hence, $V_\text{det}$ is, up to a constant multiplicative factor, the square of that volume, 
while $V_\text{det}$ is its logarithm. 
Let us write $V_\text{det}$ and $V_\text{logdet}$ as functions of the singular values of $\W$, denoted $\sigma_i(\W)$ for $1 \leq i \leq r$: 
\[
V_\text{det}(\W) 
= \frac{1}{2} \prod_{i=1}^r \sigma_i^2(\W), 
\] 
and 
\[
V_\text{logdet}(\W) 
= \sum_{i=1}^r \log(\sigma_i^2(\W) + \delta). 
\]
Both functions $V_\text{det}$ and $V_\text{logdet}$ are non-decreasing functions of the singular values of $\W$. 
This motivates us to consider 
\[
V_*(\W) = \| \W \|_* = \sum_{i=1}^r \sigma_i(\W). 
\] 
The reason is that this regularizer is also a non-decreasing function of the singular values of $\W$, and has been widely used in machine learning for several tasks such as matrix completion~\cite{recht2010guaranteed}. However, to the best of our knowledge, it has never been used in the context of VRNMF and it would be interesting to know how it compares to the standard choices $V_\text{det}$ and $V_\text{logdet}$. 
As we will see, this regularizer also performs well in practice, although not as well as  $V_\text{det}$ and $V_\text{logdet}$. 

The approach of using a volume regularization with NMF has a long history and has been considered for example in 
\cite{miao2007endmember, 
schachtner2009minimum, 
zhou2011minimum, 
fu2016robust, 
ang2018volume, 
fu2018nonnegative, 
leplatminimum}. 
The key differences among these works are in the choice of $V$. Almost all previous works have focused on the two functions $V_\text{det}$ and $V_\text{logdet}$. 
We believe it is important to design efficient algorithms for these regularizes, and to compare them on solving blind HU on highly mixed hyperspectral images in order to know which one performs better in which situations. These are the main goals of this paper. 


\subsection{Contributions}

The contribution of this work is twofold: the first part of this work is algorithm design, in which we implement and enhance the algorithms to solve VRNMF by using block  coordinate descent and optimal first-order methods.
Experimental results will show that our algorithms perform better than the current state-of-the-art volume-regularization based method from~\cite{fu2016robust}.
The second part of this work is focused on model comparisons. 
We will answer the question ``which volume function is better suited for VRNMF to tackle HU?''. 
We do so by performing extensive numerical comparisons of VRNMF with different volume functions  under various settings.
The summary of the findings are as follows:
\begin{itemize}

	\item For data that is highly separable, meaning that there exists data points close to the true endmembers, and in the presence of a few endmemebers, VRNMF with $V_\text{det}$ produces the best results in most cases. 
	
	\item For data that is less separable or in the presence of a large number of endmemebers, VRNMF with $V_\text{logdet}$ performs best in general.
	
\end{itemize}
Finally, as a proof of concept, we showcase the ability of VRNMF to produces a meaningful unmixing on hyperspectral images using real-world data. 

This work is the continuation of the conference paper~\cite{ang2018volume}. The additional contributions of this extended version are the following: 
\begin{itemize}

	\item We base our numerical experiments only on real endmembers, as opposed to randomly generated ones in~\cite{ang2018volume}.
	
	\item We use a fine grid search by bisection to tune the regularization parameter $\lambda$.
	
	\item We implement VRNMF with the nuclear norm regularizer.

\item We have improved our implementations; they are available from \url{https://angms.science/research.html}. 
	
	\item We compare our implementation of VRNMF with the state-of-the-art volume-regularization algorithm \texttt{RVolMin}~\cite{fu2016robust}. 
\end{itemize}

%
%
\subsection{Outline of the paper} \label{OrganizationNotation}

The remaining of this paper is organized as follows. 
In Section~\ref{Intro_background}, we give a more complete introduction to HU, followed by the discussion of the pure-pixel assumption also known as the separability assumption. This motivates the use of VRNMF when this assumption is violated. 
Section~\ref{Background} gives the details about the enhancement and implementation of the algorithms for VRNMF. 
Section~\ref{Experiment} presents the experiments on VRNMF and Section~\ref{Conclusion} concludes the paper.

%
%
\section{Brief review on HU}\label{Intro_background}

We now give a brief review on HU; see \cite{bioucas2013hyperspectral, ma2014signal} and the references therein for more details. 
The goal of blind HU is to study the composition and the distribution of materials in a given scene being imaged.  
A scene usually consists of a few fundamental types of materials called \textit{endmembers}, and the first goal of blind HU is to obtain the information of these endmembers from the observed hyperspectral image (HSI).  HSI are images captured by sensors over different wavelengths in the electromagnetic spectrum. 
These images form a hyperspectral data cube of size $m \times \text{col} \times \text{row}$, where $m$ is the number of spectral bands, and ``col" and ``row" are the dimensions of the images, with $n = \text{col} \times \text{row}$. 
The $m$-by-$n$ data matrix $\X$ is obtained by stacking the $m$-dimensional hyperspectral signature of the pixels as the columns of $\X$. Given the observed data matrix $\X$, the goal of blind HU is to (i)~identify the number of endmembers $r$, (ii)~obtain the spectral signature of these endmembers, (iii)~identify which pixel contains which endmember and in which proportion. 
This work does not consider problem (i) by assuming $r$ is known. In fact, model order selection is a topic of research on its own; see, e.g.,~\cite{bioucas2013hyperspectral} and the references therein. 
The focus of this work is (ii): recover the (ground truth) endmember spectral matrix, denoted by $\W^\text{true}$, from the observed data $\X$. In fact, assuming (ii) is solved, (iii) can be tackled by solving a (convex) nonnegative least squares problem~\cite{bioucas2013hyperspectral, ma2014signal}.  
The most widely used model for blind HU is the linear mixing model for which a hidden low-rank linear mixing structure for the data, namely $\X = \W^\text{true}\H^\text{true} + \N$. The data $\X$ is hence generated by $\W^\text{true}$ (the basis, where each column of $\W^\text{true}$ is the spectral signature of an endmember), weighted by the matrix $\H^\text{true}$ plus some noise $\N$. 
The matrix $\H$ in the HU literature is also called the \textit{abundance matrix} and encodes how much each endmember (columns of $\W^\text{true}$) is present in each pixel of the image: $\H_{k,j}$ is the abundance of the $k$th endmember in the $j$th pixel.  
There are two physical constraints in HU: the non-negativity constraints $\W \geq 0$ (spectral signatures are nonnegative),  and the nonnegativity $\H \geq 0$ and sum-to-one constraint $\H^\top\ones_r = \ones_n$ (the abundances in each pixel are nonnegative and sum to one). In this paper, we consider a more general model, namely using $\H^\top\ones_r \leq \ones_n$, that allows to take into account different intensities of illumination among the pixels of the image. 
Finally, NMF is the right model to perform blind HU under the linear mixing model, that is, to learn the endmember matrix $\W$ and the abundance matrix $\H$ from the data matrix $\X$. However, NMF is a difficult problem in general~\cite{vavasis2009complexity}.

\subsection{Pure-pixel assumption}

Separable NMF (SNMF) is able to solve blind HU when the data satisfies the \textit{separability} condition, which is also known as the \textit{pure-pixel assumption} in HU. 
It means the data $\X$ contains at least $r$ pure pixels where each pure pixel contains only one endmember, and there is a one-to-one correspondence between the $r$ pure pixels and the $r$ endmembers. 
Mathematically, separability changes the NMF model  (\ref{NMF}) to 
\vspace*{-1mm}
\begin{equation}
(\textrm{SNMF}): ~ \X \cong \W \underbrace{[ \I_r \, \H']\mathbf{\Pi}_n}_{\H} = [\X(:,\mathcal{A}) \; \X(:,\mathcal{A})\H']\mathbf{\Pi}_n,
\nonumber 
\vspace*{-2mm}
\end{equation}
where 
$\mathbf{\Pi}_n$ is a $n$-by-$n$ permutation matrix, 
and $\H' \in \IR^{r \times (n-r)}_+$ has its columns with $l_1$ norm smaller than one. 
In SNMF, we have that 
$\W = \X(:,\mathcal{A})$, that is, the index set $\mathcal{A}$ contains the indices of the pure pixels. 
Since all columns of $\H'$ have $l_1$ norm smaller than one, 
$\X(:,\mathcal{A})$ and the origin are the $r$ extreme points of the data cloud $\X$, that is, the convex hull of $\X(:,\mathcal{A})$ and the origin encapsulates all the other  points in $\X$.
Hence, SNMF is geometrically a vertex identification problem: given $\X$, locate the extreme points $\W = \X(:,\mathcal{A})$ which will be exactly $\W^\text{true}$ if the separability condition holds and no noise is present. 
Many algorithms exist to perform this task (referred to as pure-pixel search algorithms), e.g., vertex component analysis (VCA)~\cite{nascimento2005vertex} or the successive projection algorithm (SPA)~\cite{gillis2014fast}; see \cite{bioucas2013hyperspectral, ma2014signal} and the references therein for more algorithms and discussions. In this paper, we will compare VRNMF to SPA, as SPA is a state-of-the-art provably robust pure-pixel search algorithm. 

\begin{figure}[!t]
\centering	
\includegraphics[width=0.5\textwidth]{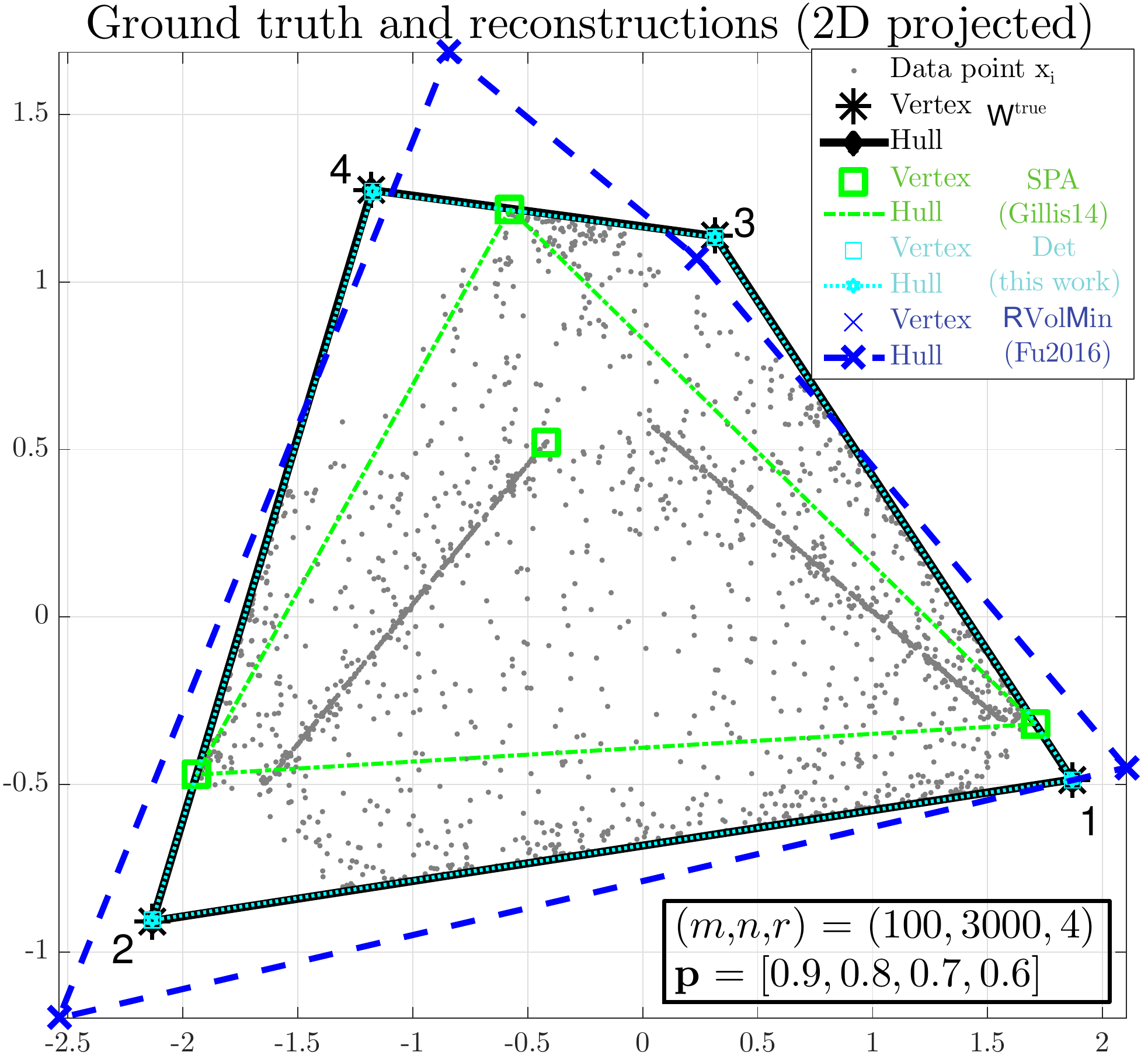}
\caption{A toy example with $(m,n,r)=(100,3000,4)$.
The plot shows the projection of data points ($\x_i$, the black dots) in two dimensions using PCA. 
This data set was generated using $\p = [0.9,0.8,0.7,0.6]$, 
meaning the maximum abundance of each ground truth vertex ($\W^\text{true}$, the black stars) in any pixel 
is at most 90\% for vertex 1, 
80\% for vertex 2, 
70\% for vertex 3, and 
60\% for vertex 4. 
}
\label{Fig_Toy}
\end{figure}

However, when the separability condition does not hold, pure-pixel search algorithms fail. 
In order to quantify how much separability is violated, we introduce the \textit{non-separability parameter} $p \in (0,1]^r$.   
First, note that separability holds if and only if each row of $\H$ contains at least one entry equal to one, that is, $\| \h^j \|_\infty = 1$ for $j \in [r]$, where $\h^j$ denotes the $j^\text{th}$ row of $\H$ and $[r] = \{1,2,\dots,r\}$.
Therefore, to break the separability condition, 
we need $\| \h^j \|\leq p_j < 1$ for some $j$. 
Figure~\ref{Fig_Toy} gives an example with $r = 4$. 
Hence, having $\| \h^j \| \leq p_j$ means that the maximum abundance of the $j$th endmember in all pixels is at most $p_j$, which sets a minimal amount of separation between the data points (the black dots in Figure~\ref{Fig_Toy}) to the vertex $\w_j$ (the black stars).  
In other words, $p_j$ controls the gap between $\{\x_i\}_{i \in [n]}$ and $\w_j$, where $\x_i$ is the $i$th data point ($i$th column of $\X$). Note that the entries of $p$ can be different in general meaning that we have \textit{asymmetric non-separability}. 
For example, in Figure~\ref{Fig_Toy}, data points are closer to vertex 1 than vertex 4. 
In the absence of pure pixels, that is, $p_j < 1$ for some $j$, it has been proved that, under mild conditions, minimizing the volume of the convex hull of the columns of $\W$ allows to recover the endmembers. These conditions require that the data points are well spread within the convex hull spanned by the columns of $\W$; see~\cite{fu2018nonnegative} for a recent survey on the topic.

Note that in Figure~\ref{Fig_Toy}, the reconstructions given by SPA~\cite{gillis2014fast} (green vertices) and \texttt{RVolMin} \cite{fu2016robust} (deep blue vertices, a  state-of-the-art minimum-volume NMF algorithm) are far away from the ground truth. 
VRNMF with $V_{\text{det}}$ (cyan vertices) produces a perfect recovery of $\W^\text{true}$.  
The next section describes how we solve the VRNMF minimization problem to obtain such results.

%
%

\section{Solving VRNMF} \label{Background}

In this section, we describe how to solve (\ref{VRNMF_pbl}).
We use the framework of block coordinate descent (BCD) by solving the subproblems on $\W$ and $\H$ separately in an alternating scheme, as done in most NMF works~\cite{gillis2014and}. 
Let us start with the subproblem for $\H$. 

%
%
\subsection{Subproblem for $\H$}\label{subsec:Solving_H}

Splitting $\H$ into columns yields $n$ independent problems: for $1 \leq j \leq n$, solve  \vspace*{-1mm}
\begin{eqnarray}
\underset{\h_j \in \Delta }{\min} \dfrac{1}{2} \| \x_j - \W \h_j \|_2^2
, ~  \Delta = \{ \h \in \IR^r_+  |  \h^\top \ones_r \leq 1 \},
\vspace*{-2mm}
\label{solving_for_H}
\end{eqnarray}
where $\h_j$ is the $j$th column of $\H$ and 
$\Delta$ is the $r$-dimensional unit simplex that encodes the non-negativity and sum-to-one constraints in (\ref{VRNMF_pbl}). 
Assuming $\text{rank}(\W) = r$ which is a standard assumption, this least squares problem over the unit simplex  is a convex problem with strongly convex objective function. 
We use the accelerated {projected gradient} (APG) method from Nesterov~\cite{nesterov2013introductory} which requires $\mathcal{O}\left(\sqrt{\kappa} \log \frac{1}{\epsilon} \right )$ iterations to reach an $\epsilon$-accurate solution, where $\kappa$ is the condition number of $\W^\top\W$ which is the Hessian of the objective function in~\eqref{solving_for_H}. 
The convergence rate of APG is optimal as no other first-order method can have a faster convergence rate~\cite{nesterov2013introductory}. 
We defer the explanation of the acceleration scheme to section \ref{subsec:logdet}.  
To compute the projection onto $\Delta$, 
we use the implementation from \cite{gillis2014successive} requiring $\mathcal{O}(r \log r)$ operations, and which uses the fact that the projection can be written as $P_{\Delta}(\h) = [ \h - l \ones_r ]_+$ where $l$ is a Lagrangian multiplier. Since $r$ is small (usually $r \leq 20$), this implementation is numerically as good as the optimal method~\cite{condat2016fast} with complexity $\mathcal{O}(r)$. 
	
In summary, we solve (\ref{solving_for_H}) using an optimal first-order method. 
Note that (\ref{solving_for_H}) is parallelizable, we can solve the $n$ problems (\ref{solving_for_H}) in parallel.

%
%


\subsection{Subproblem for $\W$ with $V_\text{det}$}\label{subsec:det}

We follow the idea from~\cite{zhou2011minimum} and perform a block coordinate descent method on the columns $\w_i$ of $\W$ ($1 \leq i \leq r$). 
We have
\begin{eqnarray}
\| \X - \W\H \|_F^2 \hspace*{-2mm} &=& \hspace*{-2mm} \| \h^i \|_2^2 \|\w_i\|_2^2  - 2\langle \X_i\h^{i \top} , \w_i \rangle + c,
\label{Lemma_bi_hals}
\\
\det(\W^\top\W) \hspace*{-2mm} &=& \hspace*{-2mm}  \gamma_i \w_i^\top \Q_i\Q_i^\top \w_i.
\label{Lemma_bi_det}
\end{eqnarray}
In (\ref{Lemma_bi_hals}), $\X_i = \X - \sum_{j\neq i} \w_j \h^j $ and $c$ is a constant independent of $\w_i$. 
In (\ref{Lemma_bi_det}), $\gamma_i = \det(\W_i^\top\W_i)$ and $\Q_i$ is the orthonormal basis of the null space of $\W_i^\top$, where $\W_i$ is $\W$ without the column $\w_i$; see \cite[Appendix]{zhou2011minimum} for more details. 
Using (\ref{Lemma_bi_hals}) and (\ref{Lemma_bi_det}), we obtain the following problem for each column of $\W$: 
\[
\min_{\w_i \geq 0}  \dfrac{1}{2} 
\w_i^\top
\big ( \| \h^i \|_2^2 \I_m + \gamma_i \Q_i\Q_i^\top \big )
\w_i 
- \langle \X_i\h^{i \top} , \w_i \rangle.
\]
Unlike the problem on $\h$, here the objective function in $\w_i$ depends on the other columns of $\W$.  
We solve this quadratic program with nonnegativity constraints 
using APG, which is faster than the standard quadratic programming algorithms used in~\cite{zhou2011minimum}.


%
%
\subsection{Subproblem for $\W$ with $V_\text{logdet}$}\label{subsec:logdet}


To solve for $V_\text{logdet}$, we use \textit{majorization minimization}, similarly as in~\cite{fu2016robust}.   
First we have
\begin{eqnarray}
(\text{Lemma 2},\cite{jose2011robust}) ~~ \log\det(\W^\top\W + \delta \I_r ) \leq  \| \W \|_\D^2 + c,
\label{TaylorBound} 
\end{eqnarray}
where $\| \W \|_\D = \| \D^{\frac{1}{2}} \W^\top \|_F $ is a weighted norm with $\D = (\Y^\top\Y + \delta \I_r)^{-1} \succ 0$ for any matrix $\Y$.
This expression (\ref{TaylorBound}) comes from performing the first-order Taylor expansion of the concave function $\log\det(.)$. 
The inequality~(\ref{TaylorBound}) holds when $\Y = \W$.
In other words, we minimize the \textit{tight convex upper bound} on the non-convex logdet function 
\begin{equation}
\underset{\W^t \geq 0}{\min}  
\Phi(\W) = \dfrac{1}{2} \langle {\W^t}^\top\W^t , \H\H ^\top \rangle - \langle \X, \W^t \H \rangle  + \dfrac{\lambda}{2} \| \W^t \|^2_{\D^t}, \label{Taylor}
\end{equation}
where $\D^t = \left( {\W^{t-1}}^\top\W^{t-1} + \delta \I_r \right)^{-1}$. 
To solve~\eqref{Taylor}, we again use APG. However, we embed the following two acceleration strategies to APG: 
\begin{itemize}
\item The adaptive restart heuristic from~\cite{o2015adaptive}.  
The APG update of $\W$ at iteration $t$ can be expressed as 
$$ 
\W^{t+1} = [\W^t - \alpha_t \nabla \Phi(\W^t + \beta_t \Delta_t) + \beta_t \Delta_t]_+, 
$$
where $\alpha$ is the step size, $\nabla \Phi$ is the gradient of the objective function and $[ \, \cdot \, ]_+$ is the projection onto the nonnegative orthant.  
Here $\beta_t \Delta_t$ is the momentum term added in Nesterov's acceleration that extrapolates $\W^t$ along the direction $\Delta_t = \W^t - \W^{t-1}$ and $\beta_t \in [0,1]$ is a parameter with $\beta_0 = 0$.
The extrapolation may increase the value of $\Phi$, and when this happens, we ``restart'' the APG scheme by reinitializing $\beta_t$, which reduces the APG back to a standard projected gradient step that decreases the value of $\Phi$.

\item The general acceleration framework for NMF algroithms from~\cite{gillis2012accelerated}. 
There are several terms independent of $\W$, in particular $\H\H^\top$ and $\X\H^\top$, in the gradient term and are the main computational cost. 
The idea has two components: 
(i)~pre-compute the terms independent of $\W$ outside the update to avoid repeated computation of the same terms, and 
(ii)~perform the update on $\W$ multiple times to reuse these precomputed terms (in most standard NMF algorithms, $\W$ is only updated once before the update of $\H$). 
\end{itemize}
Finally, due to space limit, we do not present another majorization for logdet which is based on eigenvalue inequality proposed in the conference work \cite{ang2018volume}. 
In short, that approach is another relaxation of (\ref{Taylor}) that unlocks the coupling between columns of $\W$ so that column-wise update similar to the one mentioned in section \ref{subsec:det} can be used.

%
%

\subsection{Subproblem for $\W$ with $V_*$}\label{subsec:nuclear}

The nuclear norm regularized VRNMF problem on $\W$ is
\begin{equation}
\min_{\W \geq 0} \dfrac{1}{2} \langle \W^\top\W , \H \H^\top \rangle 
- \langle \X, \W\H \rangle 
+ \lambda \| \W \|_*.
\label{nuclearVRNMF}
\end{equation}
Ignoring for now the non-negativity constraints, we solve the resulting problem by proximal gradient which updates $\W$ as follows: 
\[
\W^{t+\frac{2}{3}} = \text{prox}_{\theta \| \, \|_*} 
\left \{  \W^{t+\frac{1}{3}} \right \}, ~ \W^{t+\frac{1}{3}} := \W^t - \alpha \nabla f(\W^t),
\]
where $\alpha$ is step size and $\nabla f$ is 
the gradient of  $\frac{1}{2} \langle \W^\top\W , \H \H^\top \rangle - \langle \X, \W\H \rangle $. For the step size, we use $\alpha = 1/L$ where $L = \| \H\H^\top\|_2$ is the Lipschitz constant of the gradient. 
The proximal operator of $\| \cdot \|_*$ with step size $\theta$ has the closed-form expression given by the \textit{singular value thresholding} (SVT) operator \cite{cai2010singular}, we have
\begin{eqnarray}
\W^{t+\frac{2}{3}} = \text{SVT}_{\theta} \left \{ \W^{t+\frac{1}{3}} \right \} := \U [\Sigma - \theta \I ]_+ \V^\top, \nonumber 
\end{eqnarray}
where $\U \Sigma \V^\top$ is the SVD of $\W^{t + \frac{1}{3}}$. 
Finally, to take the non-negativity constraint into account, 
we simply put a nonnegative projection after SVT and set $\W^{t+1} = [ \W^{t + \frac{2}{3}}]_+$.

%
%
\subsection{Summary of the algorithms} 

Algorithm \ref{VRNMF_Algo} shows the general framework for solving VRNMF.
Our implementations are faster than previous works because of the use of Nesterov's acceleration~\cite{nesterov2013introductory}, 
adaptive restart~\cite{o2015adaptive} 
and the use of the acceleration strategy for NMF from~\cite{gillis2012accelerated}. 
Moreover, our update of $\W$ in VRNMF using $V_*$ is new. 

\begin{algorithm}
\caption{VRNMF}
\begin{algorithmic}[H]
\renewcommand{\algorithmicrequire}{\textbf{Input:}}
\renewcommand{\algorithmicensure}{\textbf{Output:}}
\REQUIRE $\X \in \IR^{m \times n}$, an integer $r$, $\lambda \geq 0 $
\ENSURE  $\W \in \IR_+^{m \times r}$, $\H \in \IR_+^{r \times n}$ for problem (\ref{VRNMF_pbl}).
\\ \textit{Initialize} $(\W^0,\H^0)$ by SPA \cite{gillis2014fast}
\FOR {$t = 1, 2, \dots $}
\STATE Update of $\W$
\STATE Compute and store $\H\H^\top$ and $\X\H^\top$ 
\IF {$V= V_\text{det}$} 
\STATE Update $\W$ as stated in section \ref{subsec:det}.
\ELSIF {$V= V_\text{logdet}$}
\STATE Update $\W$ as stated in section \ref{subsec:logdet}.
\STATE Update $\D = (\W^\top\W + \delta \I_r)^{-1}.$
\ELSIF {$V= V_*$}
\STATE Update $\W$ as stated in section \ref{subsec:nuclear}.
\ENDIF
\STATE Update of $\H$ by FGM \cite{gillis2014successive}.
\ENDFOR
\end{algorithmic} 
\label{VRNMF_Algo}
\end{algorithm}

We refer to the algorithm using $V_\text{det}$, $V_\text{logdet}$ and $V_*$ as \texttt{Det}, \texttt{logdet} and \texttt{Nuclear}, respectively.

\begin{figure}[!t]
\centering
\includegraphics[width=0.5\textwidth]{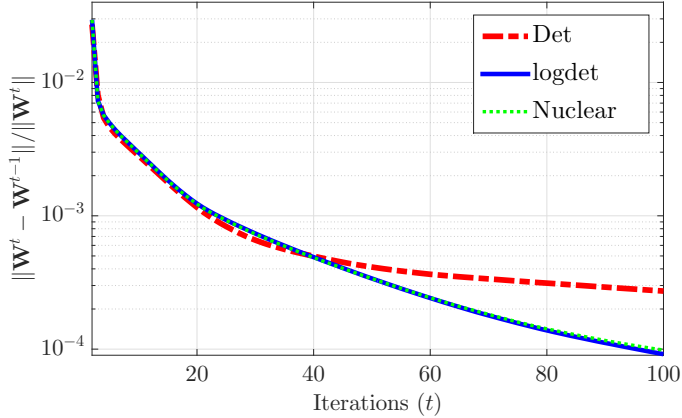}
\caption{The typical convergence curve of the algorithms \texttt{Det}, \texttt{logdet} and \texttt{Nuclear}.}
\label{Fig_convergence}
\end{figure}

We end this section by briefly mentioning the convergence of these algorithms.
VRNMF is a non-convex problem and it can be shown that the sequence $\{\W^t,\H^t\}_{t \geq 1}$ produced by 
Algorithm~\ref{VRNMF_Algo} converges to a first-order stationary point: 
For \texttt{Det}, the convergence comes from the standard result of coordinate descent~\cite{wright2015coordinate}.
For \texttt{logdet}, which involves the use of the upper bound $\Phi'$, convergence comes from the theory of  \cite{razaviyayn2013unified}. 
For \texttt{Nuclear}, convergence result of proximal gradient applies~\cite{cai2010singular}.
Figure~\ref{Fig_convergence} shows a typical convergence curve of the algorithms.

%
%
\section{Experiments}\label{Experiment}

In this section we compare the VRNMF models on solving HU.
We first describe the general experimental setting in section \ref{subsec:setting} and then report the experimental results on recovering the ground truth $\W^\text{true}$ under different asymmetric non-separability and different noise levels in the subsequent sub-sections.
Finally we present results on two real-world data sets. 
All the experiments are run with \texttt{MATLAB} (v.2015a) on a laptop with 2.4GHz CPU and 16GB RAM. The codes available from  \url{https://angms.science/research.html}. 
	
\subsection{Settings}\label{subsec:setting}

\begin{figure*}[tb]
\begin{centering}
\scalebox{0.95}{\subfloat[{\small{Samson}}\label{fig:realHyperImages_samson}]{
\begin{centering}
\includegraphics[width=0.32\columnwidth]{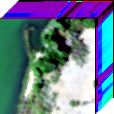}
\par\end{centering}
				
\centering{}}\subfloat[{\small{Jasper Ridge (Jasper)}}\label{fig:realHyperImages_japser}]{
\begin{centering}
\includegraphics[width=0.35\columnwidth]{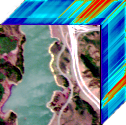}
\par\end{centering}
				
\centering{}}\subfloat[{\small{Urban}}\label{fig:realHyperImages_urban}]{
\begin{centering}
\includegraphics[width=0.48\columnwidth]{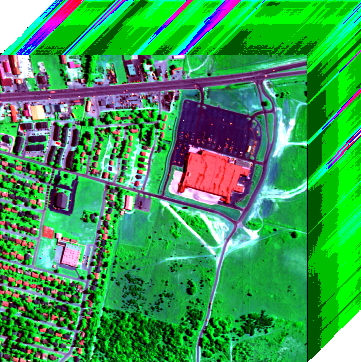}
\par\end{centering}
				
\centering{}}\subfloat[{\small{Cuprite}}\label{fig:realHyperImages_cuprite}{\small{ }}]{
\begin{centering}
\includegraphics[width=0.35\columnwidth]{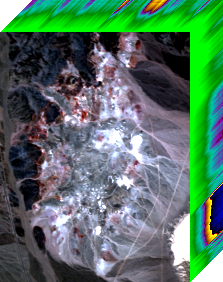}
\par\end{centering}
\centering{}}}

\centering
\includegraphics[width=170mm]{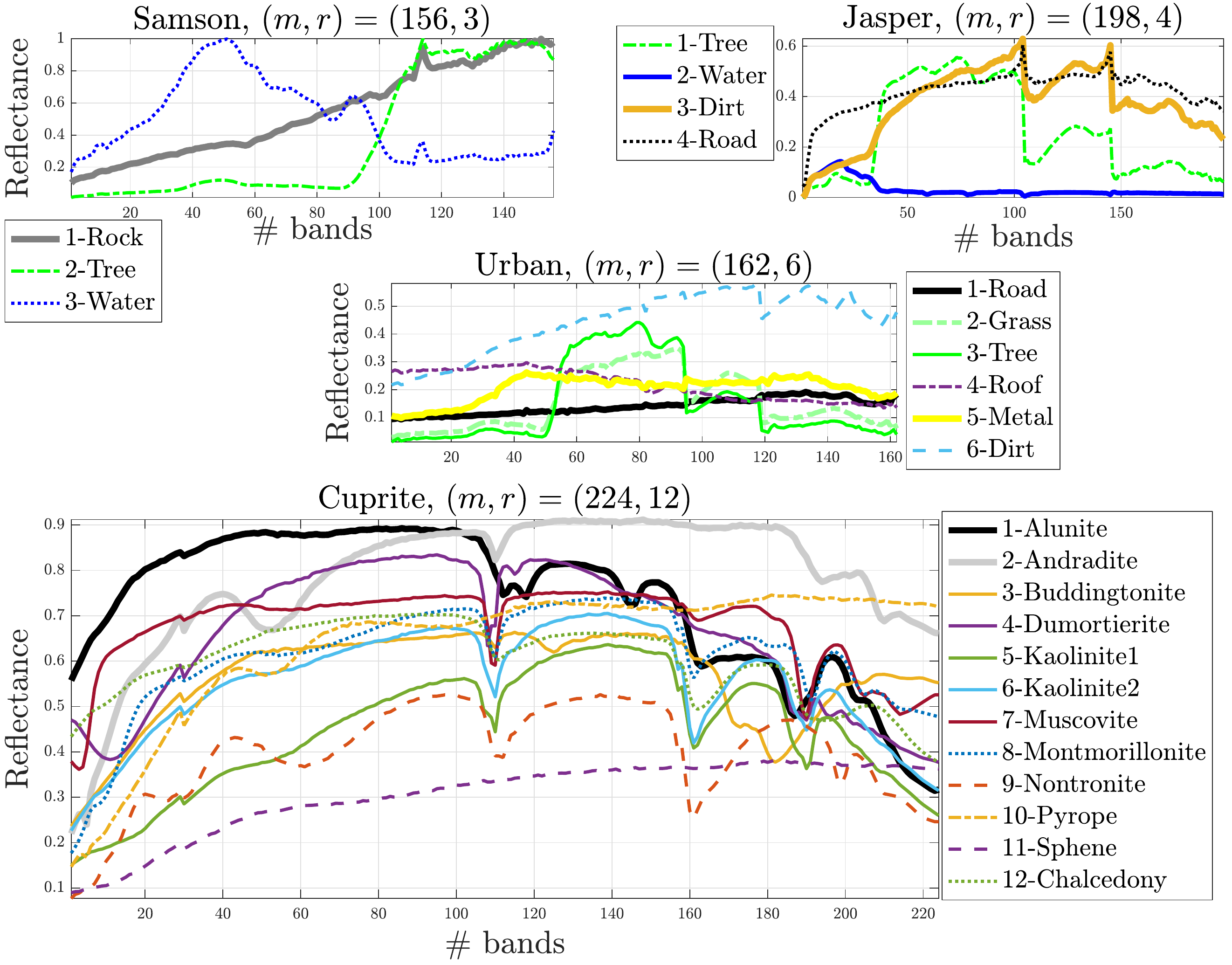}
\par\end{centering}
\centering{}
\caption{\textbf{Top row} The datasets from~\cite{zhu2014spectral}. 
\textbf{Bottom subplots} The endmembers $\W^{\text{true}}$ of the datasets.}
\label{Fig_dataset}
\end{figure*}

\paragraph{Data generation}
For synthetic experiments, we generate the observed data matrix as $\X = [\X^\text{clean} + \N]_+$ with white Gaussian noise $\N \in \IR^{m \times n}$ with zero mean and variance $\sigma \geq 0$.
 We generate $\X^{\text{clean}} = \W^\text{true} \H^\text{true}$, where $\W^\text{true}$ comes from several datasets available from \url{http://lesun.weebly.com/hyperspectral-data-set.html}~\cite{zhu2014spectral} (unlike the conference version \cite{ang2018volume} that generated $\W^{\text{true}}$ at random); see Figure~\ref{Fig_dataset}. 
We generate each column of 
$\H \in \IR_+^{r \times n}$ using the Dirichlet distribution of parameter $0.1$. 
If $\p = \ones_r$, that is, separability condition holds, we take $\H$ as $\H^\text{true}$.
Otherwise, we remove the columns of $\H$ with at least one element in the $j$th coordinate that exceeds the value $p_j$, and resample again until $\H$ satisfies the condition $\| \h^j \|_{\infty} \leq p_j$ for all $j$. 
In this paper, we will use $p_j \in (0.5, 0.99]$ for all $j$. In all the experiments, we set the number of data points $n = 1000$ and the maximum number of iterations to 300. 

For experiments on real data, we take the real data as $\X$, without any preliminary dimension reduction, and without any pre-processing to suppress noise or to remove outliers.

\paragraph{Parameters of the algorithms} 
We will compare our three proposed algorithms \texttt{Det}, \texttt{logdet} and \texttt{Nuclear} with SPA~\cite{gillis2014fast} and \texttt{RVolMin} which is a state-of-the-art volume-regularized method~\cite{fu2016robust}. 
 For \texttt{RVolMin}, we use the same data fitting term (namely, the Frobenius norm), the same number of iterations, the same parameter search scheme and the same initialization as for our algorithms. 

Given an observed data matrix $\X$, all VRNMF algorithms have two main parameters: $r$ and $\lambda$. They also require an initialization $(\W^{\text{ini}}, \H^\text{ini})$. We assume $r$ is known. We generate $\W^{\text{ini}}$ from $\X$ using SPA \cite{gillis2014fast}, and  generate $\H^\text{ini}$ using the method described in Section~\ref{subsec:Solving_H}. 
The regularization parameter $\lambda$ should usually be chosen small. In fact, a large $\lambda$ forces the vertices of the convex hull of $\W$ to be very close to each other, making $\W$ ill-conditioned and/or rank-deficient. In particular, for VRNMF with $V_*$, the SVT operator set singular values smaller than $\lambda$ to zero, so a large $\lambda$ makes $\W$ rank-deficient. 
The following describes how we tune $\lambda$. The goal here is to tune $\lambda$ so that each algorithm performs as best as possible for the considered problems. To achieve this goal, we use the ground truth $\W^{\text{true}}$.   
We set $\lambda = \tilde{\lambda} \dfrac{ f(\W^{\text{ini}},\H^\text{ini})}{ |V(\W^{\text{ini}})|}$,
where $(\W^{\text{ini}},\H^\text{ini})$ are the initial solutions obtained by SPA, and $\tilde{\lambda}$ is a tuning variable within the search interval $\mathcal{I}_0 = [10^{-6},0.5]$.
We perform grid search by bisection in a greedy way to tune $\tilde{\lambda}$: 
\begin{itemize}
	\item (step-1) Take $\tilde{\lambda}_a = 10^{-6}$, $\tilde{\lambda}_b = 0.5$ initially.
	\item (step 2) Run the algorithm with $\tilde{\lambda}_a$, $\tilde{\lambda}_b$ and $\tilde{\lambda}_c$ where $\tilde{\lambda}_c = 0.5 (\tilde{\lambda}_a + \tilde{\lambda}_b)$.
	Denote the performance of the solution $\W$ produced under a specific $\tilde{\lambda}$ as MRSA($\tilde{\lambda}$), where MRSA will be defined in the next paragraph. 
	\item (step-3) Split the search interval $\mathcal{I}_0 =[\tilde{\lambda}_a$,$\tilde{\lambda}_b]$ into two intervals $\mathcal{I}_0^{1} =[\tilde{\lambda}_a,\tilde{\lambda}_c]$ and $\mathcal{I}_0^{2} =[\tilde{\lambda}_c,\tilde{\lambda}_b]$.
	For each interval there are two MRSA values, we let the MRSA value of an interval be the sum of these values.
	\item (step-4)
	We repeat step-1 to step-3 on the interval with the lowest MRSA value.
	That is, at iteration $k$, we shrink the search interval by half by defining the new search interval $\mathcal{I}_{k+1}$ as
	$$
	\mathcal{I}_{k+1} \leftarrow \min_\mathcal{I} \left \{ \text{MRSA}(\mathcal{I}_k^1) , \text{MRSA}(\mathcal{I}_k^2) \right \}. 
	$$
	\item If a draw happens in step-4, we perform two bisections on each of the interval $\mathcal{I}_k^1$ and $\mathcal{I}_k^2$ and set the next interval to be the one with the smallest MRSA.
	\item We repeat this process (step-1 to step-4) at most 20 times, or if the improvement from one iteration to the next is negligible, namely if  
	$$ 
	\left | \text{MRSA} \left (\tilde{\lambda}_{k+1} \right ) - \text{MRSA} \left (\tilde{\lambda}_k \right  ) \right | \leq 10^{-4}. 
 	$$
\end{itemize}
Note that such greedy bisection search does not guarantee that $\tilde{\lambda}_k$ will converge to the best value, that is, the value that corresponds to the lowest MRSA.  
However, we have observed in extensive numerical experiments the effectiveness of this scheme,
which will be illustrated in \ref{subsec:LambdaSearch}.



\paragraph{Performance metric} We measure the quality of a solution $\W$ produced by VRNMF algorithms using the \textit{mean removed spectral angle} (MRSA) between $\W$ and $\W^{\text{true}}$.
MRSA between two vectors $\x , \y \in \IR^m \setminus \{\mathbf{0}_m\}$ is defined as 
\begin{equation}
\dfrac{100}{\pi} \cos^{-1} \left( \dfrac{ \langle  \x - \bar{\x}, \y - \bar{\y} \rangle }{ \| \x - \bar{\x} \|_2 \|\y - \bar{\y} \|_2}  \right) \in [0,100]. \label{mrsa}
\end{equation}
MRSA gives a better measurement than relative percentage error as it purely depends on the shapes of $\x$ and $\y$, where the effects of shifts and scaling are removed.
A low MRSA value means a good matching between $\x$ and $\y$. 
We measure the performance of algorithms by calculating the MRSA between $\{ \w_i \}_{i \in [r]}$ and $\{ \w_i^\text{true} \}_{i \in [r]}$, which is defined as the mean of MRSA between each vector pair $\{ \w_i,\w_i^\text{true}\}$, such value is within [0, 100].

%
%
\subsection{Effectiveness of the bisection search on  $\lambda$}\label{subsec:LambdaSearch}

In this section, we illustrate the effectiveness of our tuning strategy for $\lambda$.
We perform experiments on the synthetic dataset where $\W^\text{true}$ comes from the Samson dataset with $r=3$ as follows. 
We consider 6 different values of the  symmetric non-separability vector
$\p = [p_1,p_2,p_3]$ where $p_1 = p_2 = p_3$ are selected from the set $\{0.93,0.89,0.86,0.83,0.79,0.76\}$. 
We run \texttt{Det} with two parameter tuning schemes: 
(1)~the bisection search mentioned in \ref{subsec:setting},  and 
(2)~a brute-force grid search:  we run \texttt{Det} with all 100 equally spaced steps values of $\tilde{\lambda}$ in the interval $\mathcal{I}_0 = [10^{-6}, 0.5]$. 

For each value of $\p$, 10 data sets are generated randomly.   
Let us denote $\lambda^b$ the value of $\lambda$ obtained by the bisection search and $\lambda^*$ by the grid search.
Table~\ref{table_lambda_search} shows the results between comparing the bisection search and the grid search displayed as 
$\dfrac{ \text{MRSA}(\lambda^b)  
- \text{MRSA}(\lambda^*) }{\big | \text{MRSA}(\lambda^*) \big|}$ in the format of (mean$\pm$std) and the average number of iterations performed by the bisection search to reach $\lambda^b$.
From Table~\ref{table_lambda_search}, we make the following two interesting observations: 
\begin{enumerate}

\item As the purity goes down, more bisections are need to identify the best $\lambda$. This was expected since, for a high purity, the initialization (SPA) provides a good initial solution. 

\item In all 50 cases, the value of $\lambda^b$ lead to a slightly smaller MRSA value than $\lambda^*$ while requiring less computations. This illustrates the fact that the bisection is able to identify the right value of $\lambda$ and to refine the search around that value with more precision than an expensive exhaustive search. 
\end{enumerate}
	
\begin{table}[!t]
	\renewcommand{\arraystretch}{1.1}
	\caption{Comparison of MRSA values of the bisection and the grid search to obtain $\lambda^b$ and $\lambda^*$, respectively. The second column is the number of iterations needed for the bisection search to terminate. 
	}
	\label{table_lambda_search}
	\centering
	\begin{tabular}{c|cc} \hline\hline  
	$p$  & $\dfrac{\text{MRSA}(\lambda^b)  - \text{MRSA}(\lambda^*) }{\big | \text{MRSA}(\lambda^*) \big|}$  
	& $\# \text{ of iterations}$  \\ \hline\hline
	0.93 & -0.0016$\pm$0.0002                                                                                & 2$\pm$0.00       \\ \hline
	0.89 & -0.0143$\pm$0.0023                              & 6.4$\pm$0.52 \\ \hline
	0.86 & -0.0235$\pm$0.0029                               & 7.9$\pm$0.32 \\ \hline
	0.83 & -0.0429$\pm$0.0032                                   & 9$\pm$0.00        \\ \hline
	0.79 & -0.0815$\pm$0.0089                                   & 11.1$\pm$0.32  \\ \hline
	0.76 & -0.0419$\pm$0.0032                                                 & 10.5$\pm$0.53 \\ \hline\hline 
	\end{tabular}
\end{table} 


%
%
\subsection{ Comparison with MVC-NMF }\label{exp:MVC}
\begin{figure}[!t]
	\centering
	\includegraphics[width=0.5\textwidth]{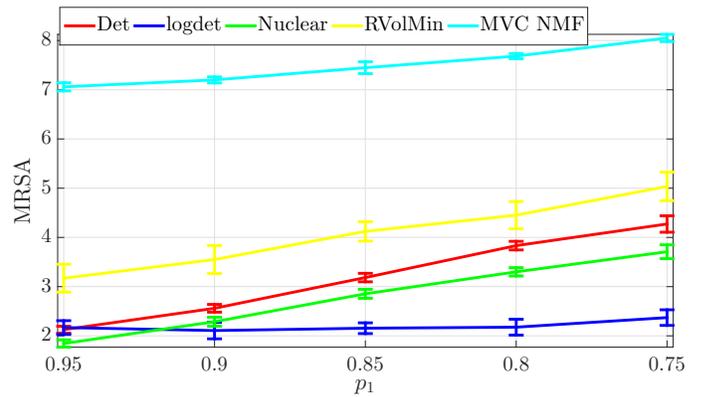}
	\caption{
		MRSA curves of \texttt{Det}, \texttt{logdet}, \texttt{Nuclear}, \texttt{RVolMin} and \texttt{MVC} across different values of $p_1$ 
	($p_2 = 0.79$, $p_3 = 0.69$) for synthetic Samson data sets. 
	} 
	\label{Fig_MVC}
\end{figure}

 Minimum volume constrained nonnegative matrix factorization (MVC-NMF) is one of the first minimum-volume NMF algorithm\footnote{The code available from \url{https://github.com/aicip/MVCNMF}}~\cite{miao2007endmember}. 
  Let us run a small experiment to show that MVC-NMF does not compete with \texttt{Det}, \texttt{logdet} and \texttt{Nuclear} and \texttt{RVolMin}~\cite{fu2016robust}.
Similarly as in the previous section, we use synthetic datasets where $\W^{\text{true}}$ comes from the Samson dataset with $r = 3$. 
However, we use more difficult scenarios using the non-separability vectors $\p =[p_1,p_2,p_3]$ where $p_1$ is selected from $\{0.95,0.9,0.85,0.8,0.75\}$, 
$p_2 = 0.79$, and $p_3 = 0.69$. 
For each value of $\p$, 10 data sets are generated randomly.   
All algorithms take the same initialization (SPA), the same number of iterations (100) and the regularization parameter $\lambda$ is tuned using the bisection search.
Figure~\ref{Fig_MVC} shows that \texttt{MVC-NMF} produces significantly worse results than the four other algorithms, hence will not be considered in the  extensive numerical experiments the subsequent comparisons. 
It is difficult to pinpoint the reasons of the poor results of MVC-NMF; we see at least two of them: 
(1)~MVC-NM uses another regularizer which is the squared volume of the convex hull of the columns of $\W$ without the origin,  
and 
(2)~it does not use optimal optimization methods (for example, it uses a fixed step size of 0.001). 

%
%


%
%
\subsection{Using the Samson dataset with $r$ = 3}\label{exp:r3}

\begin{figure}[!t]
\centering
\includegraphics[width=0.5\textwidth]{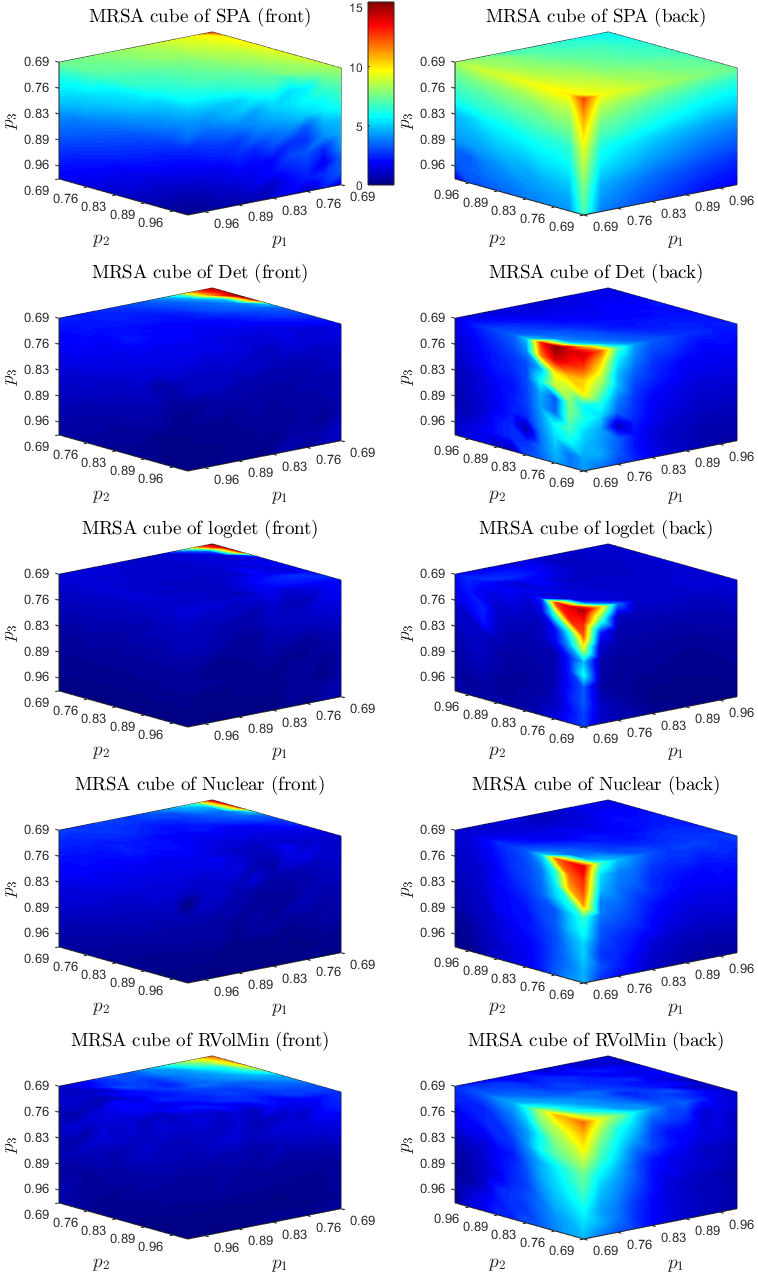}
\caption{MRSA cubes for the different NMF algorithms viewed from two different angles, for the synthetic Samson data sets with $r=3$.  
All the cubes share the same color scheme.
In general, MRSA cubes for VRNMF have a wider blue region than that of \texttt{SPA}.
However, when the vector $\p$ has a low value, VRNMF gives less accurate fitting than \texttt{SPA} (the red regions). VRNMF with \texttt{logdet} provides the best results (widest dark blue region). Best view in color.}
\label{Fig_cube}
\end{figure}

\begin{figure}[!t]
\centering
\includegraphics[width=0.5\textwidth]{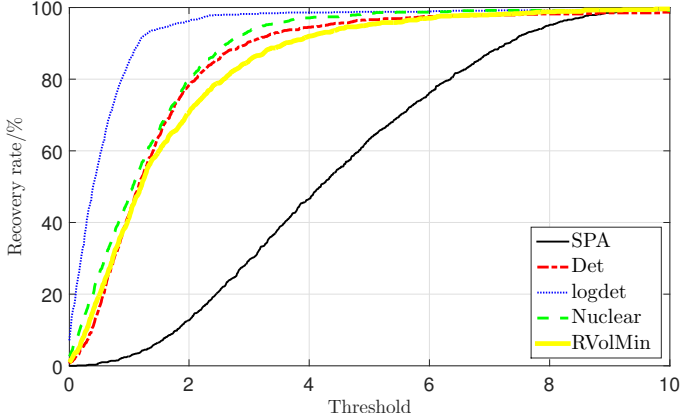}
\caption{Curves of recovery corresponding to Figure~\ref{Fig_cube} for the synthetic Samson data sets with $r=3$.}
\label{Fig_curve}
\end{figure}

Now we run more comprehensive experiments on the synthetic dataset where $\W^{\text{true}}$ comes from the Samson data set with $r = 3$ as follows. 
We perform $10\times10\times 10 = 1000$ experiments where the data in each experiment is generated using the non-separability vector $\p =[p_1,p_2,p_3]$, where $p_i$ is selected from the set $\IP = \{0.99,0.96,0.93,0.89,0.86,0.83,0.79,0.76,0.73,0.69\}$.
Figure~\ref{Fig_cube} shows the results in form of MRSA cubes for the algorithms \texttt{SPA}, \texttt{Det}, \texttt{logdet}, \texttt{Nuclear} and \texttt{RVolMin}, where each pixel in the cube is the result (in MRSA) over one trial.
We then construct the recovery curves by counting the number of cases (pixels) in the cube that the MRSA value is less than a threshold; see Figure~\ref{Fig_curve}. 
In general, the results in Figures~\ref{Fig_cube} and~\ref{Fig_curve} show that 
\begin{itemize}

	\item All VRNMF algorithms perform better than the state-of-the-art SNMF algorithm \texttt{SPA}, as the MRSA cubes of VRNMF have a wider blue region (lower MRSA) and their recovery curves in Figure~\ref{Fig_curve} dominates that of \texttt{SPA}. 
	
	\item When the data is highly non-separable (low $p_i$s), VRNMF algorithms perform worse than \texttt{SPA}. In fact, the region of the cube corresponding to highly non-separable data in \texttt{SPA} is not as red as for the VRNMF approaches. The reason is that \texttt{SPA} always extract points from the data could hence these points are never too far from the vertices. However, when the data is highly non-separable, VRNMF may generate points further away than the vertices.  
	
	

	\item Compared with \texttt{RVolMin}, \texttt{logdet} and \texttt{Nuclear} are consistently better, while \texttt{Det} performs similarly.
	 
	 \item \texttt{logdet} performs better than \texttt{Det} and \texttt{Nuclear} for highly non-separable data, as its red region is much more concentrated around the highly non-separable corner of the cube.

 
	
\end{itemize}

In terms of computational time, 
\texttt{Det} takes 2.1$\pm$0.2 seconds, 
\texttt{logdet} takes 1.2$\pm$0.1 seconds, 
\texttt{Nuclear} takes 1.3$\pm$0.2 seconds, and 
\texttt{RVolMin} takes 2.1$\pm$0.0 seconds. 
	As expected, due to the inner loop and the computation of $\Q_i$, \texttt{Det} is slower than the other algorithms.  

%
%

%
%

%
%
\subsection{Using the  Jasper dataset with $r$ = 4}
%
%
\begin{figure}[!t]
	\centering
	\includegraphics[width=0.5\textwidth]{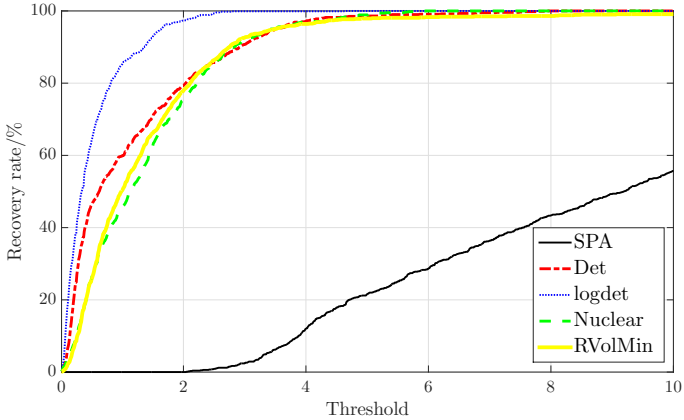}
	\caption{Curves of recovery corresponds for the dataset Jasper with $r=4$.
	Same experimental set up as the one in Figure~\ref{Fig_curve}, but where $p_4$ fixed at 0.75.}
	\label{Fig_curve_r4}
\end{figure}

In Figure~\ref{Fig_curve_r4}, we perform the same experiment as in the previous section on the dataset Jasper with $r =4$, where we fix $p_4 = 0.75$. 
Furthermore, Table~\ref{table_mrsa_Jasper} shows the result in MRSA on the dataset Jasper across three sets of predefined $\p$ values 
(highly separable with $\p_\text{high} = [0.9,0.8,0.7,0.6]$, 
less separable with $\p_\text{mid} = [0.8,0.7,0.6,0.51]$, and even less separable with $\p_\text{low} = [0.7,0.65,0.55,0.51]$)  
and three noise levels ($\sigma = 0.001, 0.005$ and $0.01$).
Each item in the table (in mean$\pm$std) is the average over 20 trails. 

We observe the following 
\begin{itemize}

\item Figure~\ref{Fig_curve_r4} show that VRNMF algorithms are competitive with the state-of-the-art minimum volume based method \texttt{RVolMin}, where \texttt{Det} performs slightly better than \texttt{RVolMin} while \texttt{logdet} is significantly better.

\item Tables \ref{table_mrsa_Jasper} shows that all the methods improve the fitting accuracy of SPA, with \texttt{Det} has the best performance in all cases and \texttt{RVolMin} has the worst performances in all cases.


\end{itemize}

In terms of computational time, 
 \texttt{Det} takes 6.62$\pm$0.5 seconds, 
 \texttt{logdet} takes 2.22$\pm$0.2 seconds, 
 \texttt{Nuclear} takes 1.45$\pm$0.0 seconds, and 
 \texttt{RVolMin} takes 2.7$\pm$0.0 seconds. 
 Hence, for the same reasons as in the previous experiment, \texttt{Det} is slower.

\begin{table}[!t]
\renewcommand{\arraystretch}{1.1}
\caption{MRSA values for the  Jasper dataset with $r = 4$}
\label{table_mrsa_Jasper}
\centering
\begin{tabular}{c|cccccc}
\multicolumn{4}{c}{Across different $\p$ ($\sigma = 0.001$)}					 \\ \hline\hline		
Method & $\p_\text{high}$      & $\p_\text{mid}$       & $\p_\text{low}$         \\ \hline\hline
SPA    & 5.40$\pm$0.60         & 12.62$\pm$0.18	       & 20.76$\pm$0.23 		 \\ \hline
Det    & \textbf{0.41$\pm$0.08}& \textbf{0.40$\pm$0.06}& \textbf{10.99$\pm$1.68} \\ \hline
logdet & 0.48$\pm$0.54 	       & 3.03$\pm$0.28		   & 12.57$\pm$1.49 		 \\ \hline
Nuclear& 0.64$\pm$0.07		   & 2.12$\pm$0.13	       & 19.90$\pm$1.40          \\ \hline
RVolMin& 1.04$\pm$0.38  	   & 4.99$\pm$0.22		   & 13.67$\pm$2.99 		 \\ \hline\hline 
\end{tabular}

\vspace*{2mm}
\begin{tabular}{c|ccccc}
\multicolumn{4}{c}{Across different noise levels ($\p = \p_\text{high}$)}         \\ \hline\hline
Method  & $\sigma = 0.001$       &  $\sigma = 0.01$      & $\sigma = 0.05$        \\ \hline\hline
SPA     & 5.40$\pm$0.60          & 7.29$\pm$0.51	     &	24.59$\pm$1.43        \\ \hline
Det     & \textbf{0.41$\pm$0.08} & \textbf{0.74$\pm$0.06}&	\textbf{4.90$\pm$3.27}\\ \hline
Taylor  & 0.48$\pm$0.54 	     & 1.40$\pm$0.08		 &	9.00$\pm$3.88         \\ \hline
Nuclear & 0.64$\pm$0.07		     & 1.23$\pm$0.06	     &	6.78$\pm$5.78         \\ \hline
RVolMin & 1.04$\pm$0.38  		 & 2.31$\pm$0.28		 &	10.23$\pm$2.05        \\ \hline\hline
\end{tabular}
\end{table}

%
%

\subsection{Using the Urban dataset with $r$ = 6} 

Table~\ref{table_mrsa_urban} shows the result of the same experiment as in the previous section performed with the Urban dataset. 
Here $\p_\text{high} = [0.9,0.75,0.7,0.65,0.8,0.85]$, $\p_\text{mid} = [0.8,0.7,0.65,0.6,0.75,0.8]$, and $\p_\text{low} = [0.7,0.6,0.55,0.51,0.65,0.7]$.
The noise levels are $\sigma = 0.001, 0.005$ and $0.01$.

The results show that 
\begin{itemize}
\item All the methods improve the fitting accuracy of SPA.

\item \texttt{Det} has the best performance in most cases, with \texttt{logdet} as the first-runner up.

\item \texttt{Det} performs well for $\p_\text{high}$ and $\p_\text{mid}$. 
	With $\p_\text{low}$, it is not as good as \texttt{logdet}.
	
\item \texttt{Nuclear} and \texttt{RVolMin} have the worst performances in all cases.

\item The computational times are: \texttt{Det} 7.5$\pm$0.2, \texttt{logdet} 1.8$\pm$0.1, \texttt{Nuclear} 1.5$\pm$0.0 and \texttt{RVolMin} 2.3$\pm$0.0, respectively.
\end{itemize}
	
\begin{table}[!t]
\renewcommand{\arraystretch}{1.1}
\caption{MRSA values for the Urban dataset with $r = 6$}
\label{table_mrsa_urban}
\centering

\begin{tabular}{c|cccccc}
\multicolumn{4}{c}{Across different $\p$ ($\sigma = 0.001$)}				\\ \hline\hline		
Method  & $\p_\text{high}$       & $\p_\text{mid}$       &  $\p_\text{low}$ \\ \hline\hline
SPA     & 7.83$\pm$0.93          & 10.32$\pm$1.53	     &	16.22$\pm$2.00  \\ \hline
Det     & \textbf{0.54$\pm$0.11} & \textbf{2.45$\pm$1.25}&	10.08$\pm$5.71  \\ \hline
logdet  & 1.27$\pm$0.68 	     & 3.09$\pm$2.04		 &	\textbf{8.78$\pm$2.58}   \\ \hline
Nuclear & 3.79$\pm$0.62		     & 6.39$\pm$1.54	     &	13.48$\pm$4.53  \\ \hline
RVolMin & 3.03$\pm$0.90  		 & 5.64$\pm$1.29		 &	13.05$\pm$4.28  \\ \hline\hline 
\end{tabular}

\vspace*{2mm}
\begin{tabular}{c|ccccc}
\multicolumn{4}{c}{Across different noise levels ($\p = \p_\text{high}$)}	\\ \hline\hline
Method  & $\sigma = 0.001$       &  $\sigma = 0.005$      & $\sigma = 0.01$ \\ \hline\hline
SPA     & 7.83$\pm$0.93          & 8.56$\pm$0.86	      &	15.95$\pm$3.57	\\ \hline
Det     & \textbf{0.54$\pm$0.11} & \textbf{1.25$\pm$0.48} &	\textbf{6.85$\pm$3.47} \\ \hline
logdet  & 1.27$\pm$0.68 	     & 4.41$\pm$0.93		  &	9.85$\pm$3.36   \\ \hline
Nuclear & 3.79$\pm$0.62		     & 4.58$\pm$0.95	      &	11.75$\pm$3.35	\\ \hline
RVolMin & 3.03$\pm$0.9  		 & 4.95$\pm$1.01		  &	15.32$\pm$9.35  \\ \hline\hline
\end{tabular}
\end{table}

%
%
\subsection{On the Cuprite dataset with $r$ = 12}

Table~\ref{table_mrsa_different_p_r12} shows the result on the same experiments conducted with the Cuprite dataset.
Here we concatenate the $\p$ vector used in Urban two times to form the $\p$ vector with length $12$ for the data Cuprite. 
For example, for $\p_\text{high}^{\text{Cuprite}}$, we use $\p_\text{high}^{\text{Cuprite}} = [\p_\text{high}^{\text{Urban}} ~ \p_\text{high}^{\text{Urban}}]$. 

The result from the two tables show that 
\begin{itemize}
\item All the methods improve the fitting accuracy of SPA.

\item \texttt{logdet} has the best performance in all situations, with \texttt{Det} and \texttt{Nuclear} as the first-runner ups.

\item \texttt{RVolMin} has the worst performances in most cases.

\end{itemize}

In terms of computational time, 
\texttt{Det} takes 28.1$\pm$0.9 seconds, 
\texttt{logdet} takes 3.8$\pm$1.1 seconds, 
\texttt{Nuclear} takes 3.4$\pm$0.1 seconds, and 
\texttt{RVolMin} takes 3.2$\pm$0.0 seconds. 
As expected, when $r$ increase, \texttt{Det} takes much more time and it is not recommended for large $r$. 

\begin{table}[!t]
\renewcommand{\arraystretch}{1.1}
\caption{MRSA values for the Cuprite dataset with $r = 12$}
\label{table_mrsa_different_p_r12}
\centering
\begin{tabular}{c|cccccc}
\multicolumn{4}{c}{Across different $\p$ ($\sigma = 0.001$)}				\\ \hline\hline		
Method  & $\p_\text{high}$       & $\p_\text{mid}$       &  $\p_\text{low}$ \\ \hline\hline
SPA     & 6.59$\pm$0.98          & 8.87$\pm$1.42	     & 11.53$\pm$1.39	\\ \hline
Det     & 2.59$\pm$0.74          & 4.40$\pm$1.51         & 9.71$\pm$2.43	\\ \hline
logdet  & \textbf{2.51$\pm$0.59} & \textbf{3.85$\pm$0.96}& \textbf{8.41$\pm$2.04}\\ \hline
Nuclear & 2.55$\pm$0.70		     & 4.33$\pm$1.56	     & 10.07$\pm$2.73   \\ \hline
RVolMin & 3.72$\pm$2.41  		 & 5.39$\pm$2.29		 & 10.82$\pm$3.31   \\ \hline\hline
\end{tabular}

\vspace*{2mm}
\begin{tabular}{c|ccccc}
\multicolumn{4}{c}{Across different noise levels ($\p = \p_\text{high}$)}   \\ \hline\hline
Method  & $\sigma = 0.001$       &  $\sigma = 0.005$      & $\sigma = 0.01$ \\ \hline\hline
SPA     & 6.59$\pm$0.98          & 7.24$\pm$1.07	     & 11.53$\pm$1.39	\\ \hline
Det     & 2.59$\pm$0.74          & 4.34$\pm$1.74         & 9.71$\pm$2.43	\\ \hline
logdet  & \textbf{2.51$\pm$0.59} & \textbf{4.01$\pm$1.22}& \textbf{8.41$\pm$2.04}\\ \hline
Nuclear & 2.55$\pm$0.70		     & 4.38$\pm$1.72	     & 10.07$\pm$2.73   \\ \hline
RVolMin & 3.72$\pm$2.41  		 & 4.66$\pm$1.73		 & 10.82$\pm$3.31   \\ \hline\hline
\end{tabular}
\end{table}	

In general, the results show that the VRNMF algorithms with \texttt{Det} and \texttt{logdet} performs very well, better than \texttt{RVolMin} and \texttt{Nuclear} in terms of fitting accuracy. 
As \texttt{RVolMin} consistently produce inferior results, we do not include it in the subsequent sections.

\subsection{On image segmentation on real data}

\begin{figure}[!t]
\centering
\subfloat[Samson - spectral signatures and abundance maps of each  component]{\includegraphics[width=90mm]{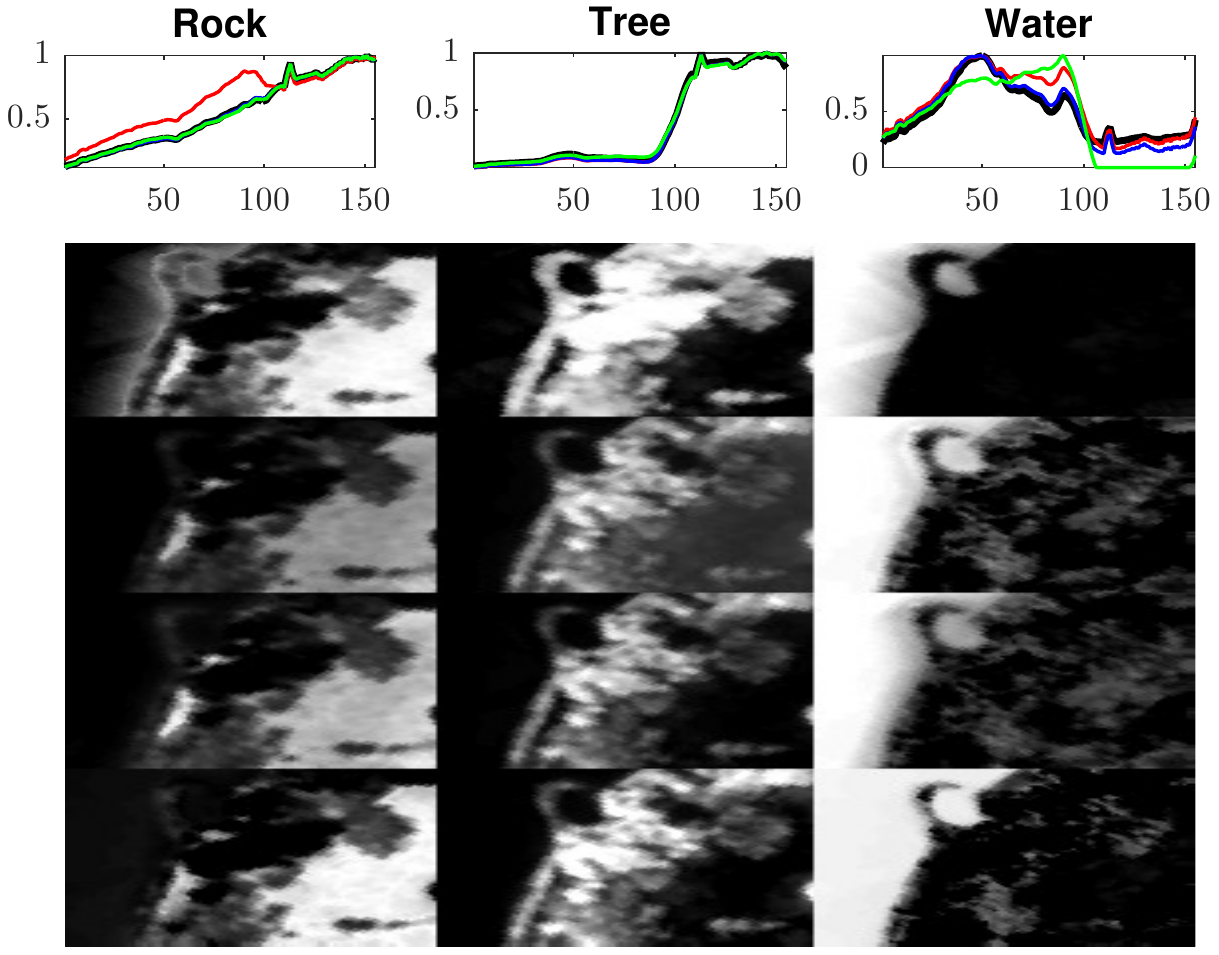}%
\label{Fig_real_Samson_abmap}} 
\\
\subfloat[Jasper - spectral signatures and abundance maps of each component]{\includegraphics[width=90mm]{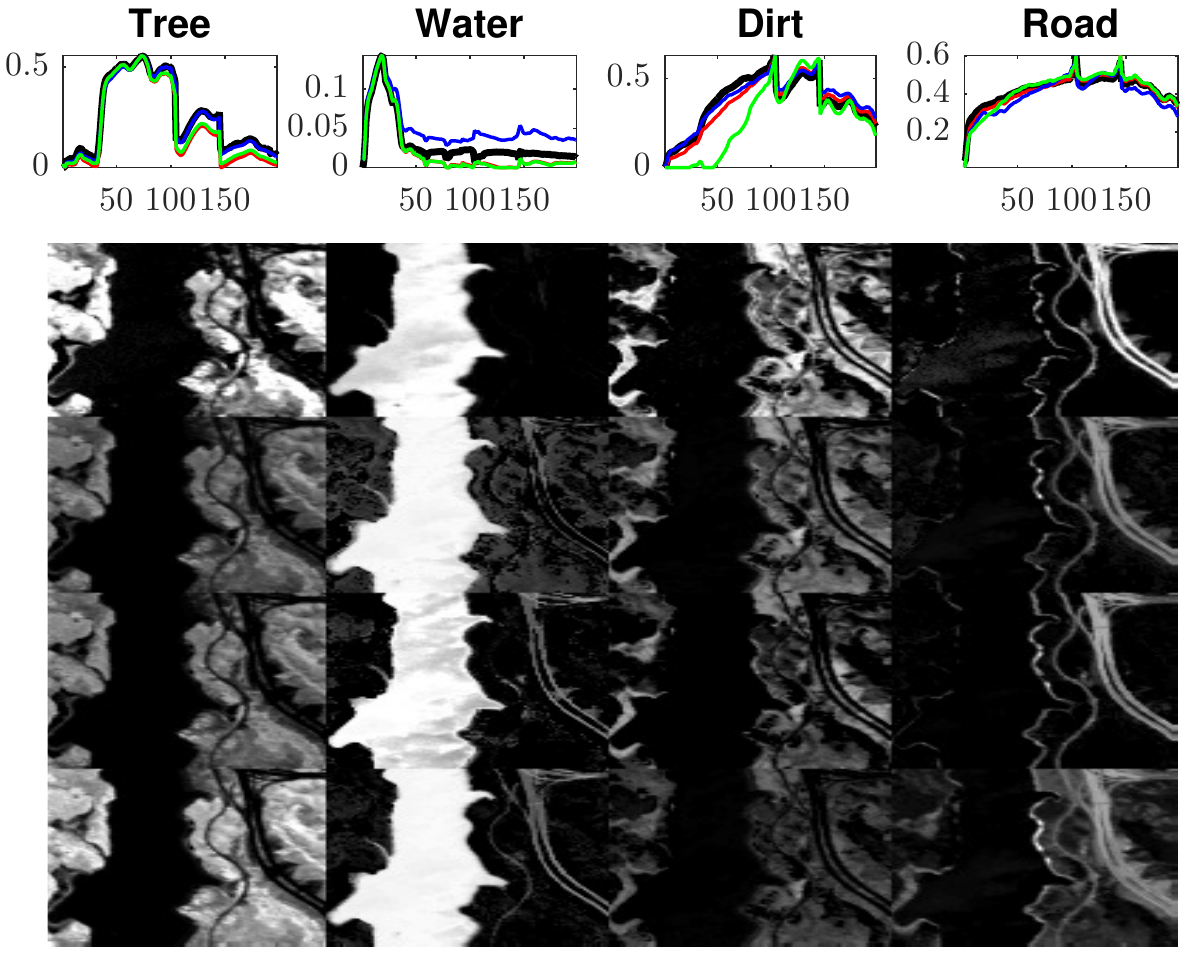}%
\label{Fig_real_Jasper_abmap}}

\caption{The decomposition of the datasets Samson and Jasper. 
In the sub-figures (a) and (b), the upper subplots are the identified endmembers $\W$: $\W^\text{ref}$ (black), $\W^\text{Det}$ (red), $\W^\text{logdet}$ (blue) and $\W^\text{Nuclear}$ (green); and the lower subplots are the corresponding abundances matrices $\H$; from the first row to the fourth row: $\H^\text{ref}$ (provided by \cite{zhu2014spectral}), $\H^\text{Det}$, $\H^\text{logdet}$ and $\H^\text{Nuclear}$. 
}
\label{Fig_real}
\end{figure}

In this section we run the algorithms \texttt{Det}, \texttt{logdet} and \texttt{Nuclear} on real HU data Samson and Jasper. 
As stated in section \ref{subsec:setting}, no pre-processing is performed on the raw data and we use the raw data directly.
The following states the specifications of the datasets.
For the dataset Samson, $(m,n,r) = (156,95^2,3)$.
For the dataset Jasper, $(m,n,r) = (198,100^2,4)$. 
We have tuned $\lambda$ in the same way as for the synthetic datasets using the endmembers $\W^\text{Ref}$ from~\cite{zhu2014spectral}. 
Figure~\ref{Fig_real} shows the decompositions. 
In the same figure, we also show the references provided by \cite{zhu2014spectral}, and we list the numerical results in Table~\ref{table_real}. 
In Table~\ref{table_real}, the MRSA is calculated with respect to the reference $\W^\text{Ref}$ of \cite{zhu2014spectral}. 
The results show that all three VRNMF algorithms produce meaningful decomposition. 

\noindent \textit{Remark}
The reference \cite{zhu2014spectral} used a sparsity regularization on $\H$ and hence the abundance map of \cite{zhu2014spectral} looks cleaner. 
It is out of the scope of this paper to consider a sparsity regularization on $\H$ as we are focusing on the volume regularization on $\W$. This is a direction for further research. 


\begin{figure}[!t]
	\centering
	\subfloat[Samson]{\includegraphics[width=0.5\textwidth]{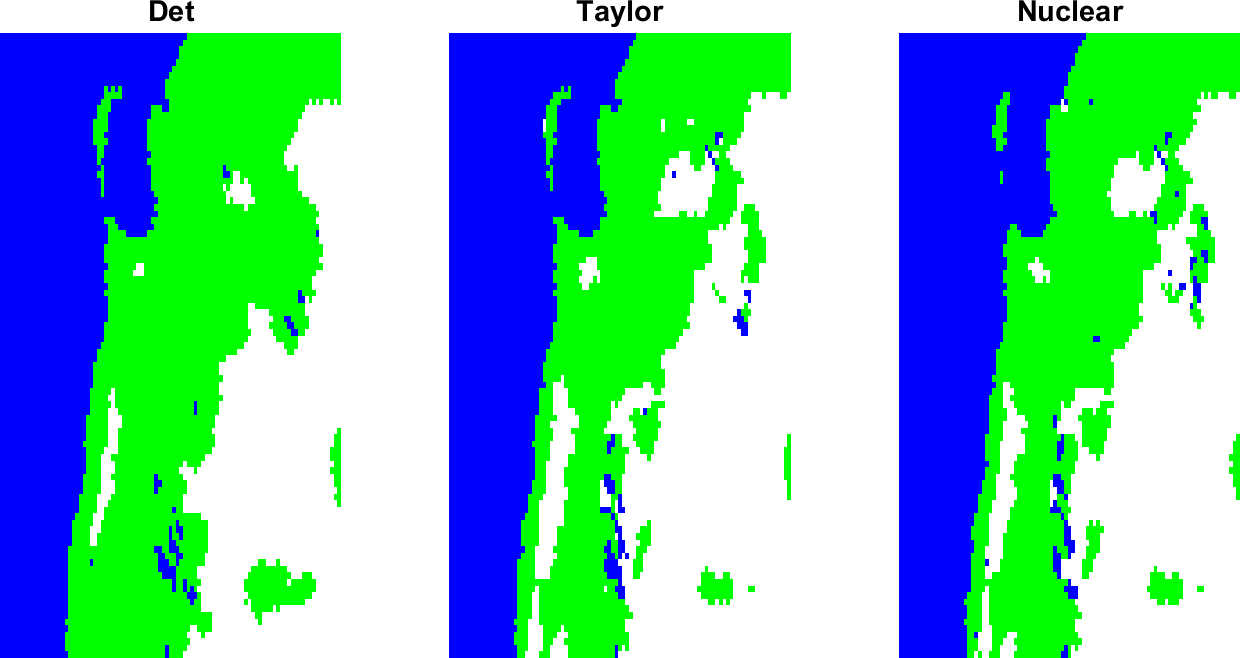}%
		\label{Fig_real_Samson}}
	\\
	\subfloat[Jasper]{\includegraphics[width=0.5\textwidth]{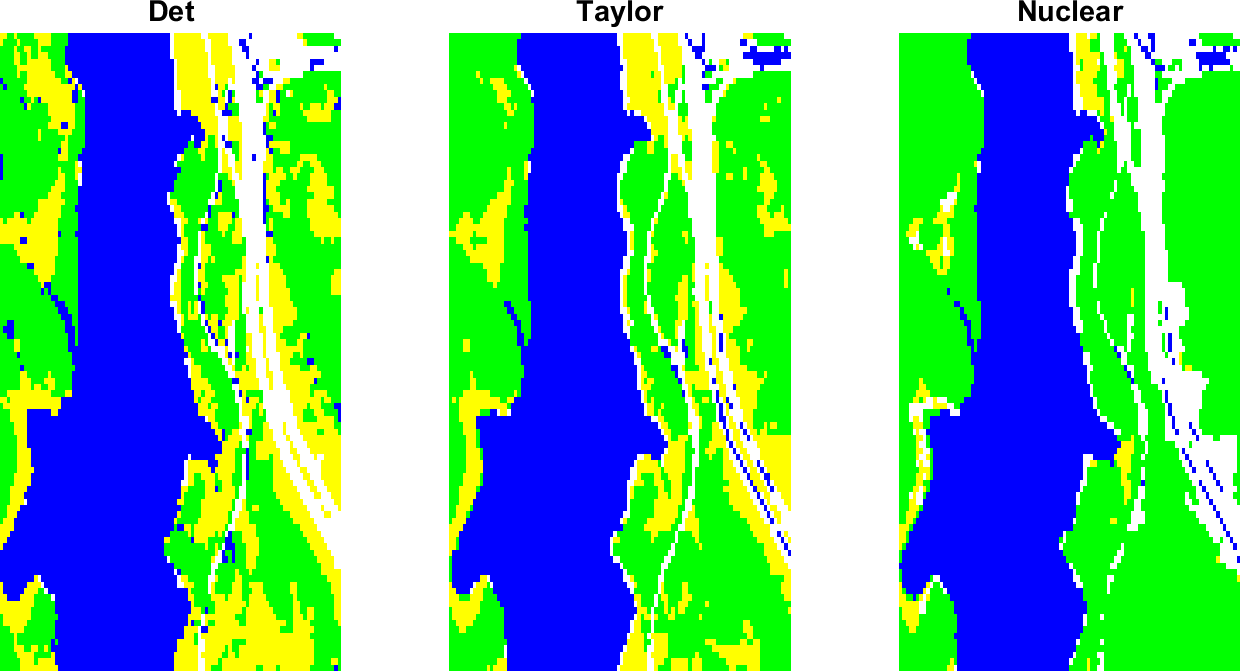}%
		\label{Fig_real_Jasper}}
	\caption{The distribution maps of the material in the whole scene obtained from the VRNMF algorithms.
		In (a), the colour code is: white -- rock, blue -- water, and green -- tree.
		In (b), the colour code is: green -- tree, blue -- water, yellow -- dirt and white -- road.
		The map is constructed by taking the largest component of $\h_j$ as the material of that pixel.
		As the factorization rank $r$ here is small and these datasets are highly separable, \texttt{Det} produces the best segmentation results. 
		For example, the road in Jasper produced by \texttt{Det} is better identified than for the other two. 
	}
	\label{Fig_real_map}
\end{figure}

\begin{table}[!t]
\renewcommand{\arraystretch}{1.1}
\caption{Numerical results on Samson and Jasper datasets.}
\label{table_real}
\centering
\begin{tabular}{c|cccc} 
\multicolumn{4}{c}{Samson} \\ \hline\hline
Method  & Time (s.) & MRSA & $\dfrac{\|\X - \W\H \|_F}{\| \X \|_F}$     \\ \hline\hline
Det     & 8.53    	& 7.13 & 2.86\%  \\ \hline
logdet  & \textbf{6.22}    	& \textbf{2.58} & \textbf{2.69}\%  \\ \hline
Nuclear & 8.94   	& 6.99 & 7.13\%  \\ \hline\hline
\multicolumn{4}{c}{Jasper} \\ \hline\hline	
Det     & 15.13    	& \textbf{5.51} & \textbf{4.14}\%  \\ \hline
logdet  & 12.67    	& 6.03 & 6.09\%  \\ \hline
Nuclear & \textbf{12.43}   	& 8.96 & 4.51\%  \\ \hline\hline
\end{tabular}
\end{table}	


\section{Conclusion}\label{Conclusion}

In this paper, NMF models with different volume regularizations (VRNMF) were investigated. We have developed highly efficient algorithms for these VRNMF models.
The VRNMF algorithms are shown to be able to outperform both the state-of-the-art separable NMF algorithm SPA, and the volume-based methods \texttt{MVC-NMF}~\cite{miao2007endmember} and \texttt{RVolMin}~\cite{fu2016robust}. 
Furthermore, extensive experimental results using real hyperspectral data showed that, when the data has a small rank $r$ ($r \leq 4$) and is highly separable, the volume regularizer based on the determinant (\texttt{Det}) provides the best results, although the 
the regularizer based on the logarithm of the determinant (\texttt{logdet}) provides almost as good decompositions. 
When the rank increases and/or the data becomes less separable, \texttt{logdet} performs best in most cases, while being computationally faster than \texttt{Det}. 
Therefore, in practice, we recommend the use of \texttt{logdet}. 


\bibliographystyle{Paper_IEEE_VRNMF_BibStyle}
\bibliography{Paper_IEEE_VRNMF_Bib}
\end{document}